\documentclass[12pt, letterpaper, ]{article}

\usepackage[utf8]{inputenc}
\usepackage{amssymb}
\usepackage{multicol}
\usepackage[margin=2.15cm]{geometry}
\usepackage{setspace}
\usepackage{graphicx}
\usepackage{url}
\usepackage{color}
\usepackage[hidelinks]{hyperref}
\graphicspath{ {images/} }
\usepackage{tikz}
\usepackage{amsmath}

\begin{document}
	
	\title{Hankel's principle as an anti-Kantian program \bigskip \bigskip}
	
	\author{\bigskip Iulian D. Toader\footnote{Institute Vienna Circle, University of Vienna, iulian.danut.toader@univie.ac.at}}
	
	\date{\bigskip (forthcoming in \textit{Archive for History of Exact Sciences})}
	
	\maketitle
	
	\doublespacing
	
	\begin{quote} 
		
		\textbf{Abstract:} Hankel used his \textit{principle of the permanence of formal laws} (PFL) as a guide for the extension of number systems and as a necessary condition for the legitimacy of their formal theories. He acknowledged that these applications have important limitations, evidenced by the extension to hypercomplex numbers and by what he saw as the unavoidable inconsistency of a formal theory of irrational numbers. Yet, intriguingly enough, he remained fully committed to the PFL. I argue that this was due to his understanding it as an expression of a conservative strategy, inherited from Peacock and Hamilton, which permits the revision of the basic laws of arithmetic if there are reasons for revision that are found, upon deliberation, to outweigh the reasons for their preservation. Then I discuss criticisms by Schubert and Pringsheim, who reformulated the PFL to align it with their own anti-revisionary conservative strategy, at the cost of relinquishing parts of modern mathematics. I conclude by emphasizing the deep philosophical difference between these kinds of conservatism in mathematics.
		
	\end{quote}
	
	
	\newpage 
	
	\section{Introduction}
	
	My goal in this paper is to offer an analysis of Hermann Hankel's applications of the principle of the permanence of formal laws (PFL) which he articulated in his \textit{Theorie der complexen Zahlensysteme} (Hankel 1867), and to determine the extent to which it can be interpreted as an anti-Kantian program in the philosophy of mathematics. Hankel used the PFL as a guide for the extension of number systems, but also more innovatively, as a guide for the introduction of their formal theories. The latter application led him to advocate a formal treatment of all mathematics, detached from any intuition of extensive magnitudes. This is an important part of his anti-Kantian program, though hardly unique to Hankel. Other mathematicians had pushed back against Kant's view that intuition is indispensable to mathematics. Bolzano and Frege, for example, objected in various ways to the use of extra-conceptual resources in proofs (Bolzano 1810, Frege 1884). However, another, arguably more significant part of this program is revealed by his view of the irrational and hypercomplex numbers as exceptions to the PFL. Intriguingly enough, Hankel did not see them as counterexamples that invalidate it or weaken its normative force. I want to argue that this was due to his understanding the PFL as an expression of a conservative strategy that allows, under certain conditions, the revision of the basic laws of arithmetic. 
	
	More specifically, on the interpretation that I propose, Hankel conceived of the PFL as an expression of \textit{deliberative conservatism}, a view that advocates the preservation of laws to the furthest extent possible, ``\textit{soweit als möglich}'' (Hankel 1867, 102). This view had been exemplified in the works of Peacock and Hamilton (Peacock 1833, Hamilton 1853; see Toader 2026), and then embraced by Hankel, who considered violations of the basic laws of arithmetic, such as distributivity or commutativity, for alternative systems of hypercomplex numbers. It also enabled him to remain fully committed to the PFL despite rejecting the formal treatment of the irrationals as impossible. I think that deliberative conservatism is an important part of Hankel's anti-Kantian program, since it is directly opposed to the kind of conservatism associated to a Kantian view of mathematics, according to which arithmetic is unconditionally valid and immune to revision (Kant 1781, Schultz 1791). This contrast proves helpful for understanding certain criticisms that later mathematicians, like Hermann Schubert and Alfred Pringsheim, leveled against Hankel.
	
	The structure of the paper is as follows. Section 2 presents Hankel's principle and distinguishes its old and more familiar applications from the new ones that he introduced. Then it focuses on the irrationals, it explains why he thought that they raised difficulties for the PFL and describes these difficulties in some detail. Section 3 clarifies Hankel's view on the hypercomplexes, and argues that he understood the PFL as an expression of deliberative conservatism, which confers the principle normative force despite its acknowledged limitations. Section 4 discusses Schubert's and Pringsheim's reformulations of the PFL and the kind of conservatism that they advocated. Section 5 concludes the paper by emphasizing the deep philosophical difference between these kinds of conservatism, and asks further questions about the relationship between formalization and conservatism in mathematics more generally.

	\section{Hankel's Principle}
	
	Hermann Hankel (14 February 1839 -- 29 August 1873) studied under, among others, Möbius, Riemann, Weierstrass, and Kronecker at the Universities of Leipzig, Göttingen, and Berlin. As a student, Hankel distinguished himself early on by winning an academy prize with a paper in which he answered a question by Dirichlet (Hankel 1861). The question asked for a derivation, in the Lagrangian formulation of hydro-dynamics, of the laws of vorticity for incompressible fluids, which had been previously proved by Helmholtz only in its Eulerian formulation (see Frisch and Villone 2014). As a young Professor of Mathematics, Hankel then taught at the University of Erlangen--Nuremberg and the University of Tübingen, where his research was mainly focused on hypercomplex numbers (Hankel 1867), the theory of functions with complex variables (Hankel 1870), as well as the history of mathematics (Hankel 1869, Hankel 1874). 
	
	Hankel's most significant contribution to the philosophy of mathematics is the adoption and application of the principle of the permanence of equivalence forms, which had been introduced several decades before by the British mathematician George Peacock (Peacock 1830, 1833, 1842, 1845). Hankel reformulated this principle and renamed it as \textit{the principle of the permanence of formal laws} (Hankel 1867, vii), which was surely meant to indicate a restriction of its application specifically to laws, and thus a deviation from Peacock, who had programmatically included theorems within the scope of his principle. The PFL was stated as follows:
	
	\begin{quote}
		\singlespacing
		
		If two forms expressed in the general signs of arithmetica universalis are equal to one another, then they shall also remain equal to each other when the signs cease to signify simple magnitudes and as a result the operations, too, acquire any content whatsoever. (1867, 11)
		
	\end{quote}
	
	\noindent An example of forms that belong to universal arithmetic, where the quantities denoted by its number signs are the positive integers, is given by the law of distributivity of multiplication over addition: $a(b+c) = ab+ac$. Another, by the laws of commutativity: $a+b = b+a$ and $ab = ba$. According to the PFL, such basic laws should be preserved in all extensions beyond the positive integers. Hankel properly acknowledged Peacock as the source of the PFL, but pointed out that his own application of it was ``new and independent.'' (1867, vii) To understand what he thought was new, it would be useful to distinguish two different applications: according to one, the PFL says that a law which holds for positive integers should hold for all other \textit{actual} numbers; according to the other -- and this was entirely new -- the PFL says that a law which holds for actual numbers should hold for \textit{formal} numbers as well. Actual (``\textit{actuelle}'') numbers are understood as expressions of relations between magnitudes (Förstemann 1817; see Schubring 2005, Salazar 2023), or more exactly relations between objects of thought that can be presented intuitively as magnitudes.  In contrast, Hankel's formal numbers are understood as expressions of relations between objects of thought that cannot be intuitively presented as magnitudes. For him, however, this was not a rigid distinction: numbers that are not known at first to be either actual or formal ``may be called potential insofar as they can be made actual, or intellectual or mental insofar as they are at first only to be thought and not intuited, or simply formal insofar as they express only a certain formal relation.'' (1867, 7sq). The new application of the PFL led Hankel to advocate a formal treatment of all numbers, detached from any intuition of magnitudes:
	
	\begin{quote} \singlespacing
		
		Through this principle, it became possible to replace the initial concept of a number, as the expression of actual relations between objects and their operations, with the more general concept of formal operations moving solely within the domain of logical thought, and of numbers emerging from the mental combination of objects -- numbers that are initially contentless, purely abstract forms of thinking. (1867, 47)
		
	\end{quote}
	
	\noindent In fact, the formal treatment he advocated, elevated to a doctrine, was far more ambitious, and underlies the first part of Hankel’s anti-Kantian program: 
	
	\begin{quote} \singlespacing
		The doctrine of forms does not have merely the narrow aim of explaining and rigorously deducing the ordinary arithmetica universalis with its whole, fractional, irrational, negative, and imaginary magnitudes; rather, with its principle of permanence, it also proves to be eminently fruitful for the entire organism of mathematics. (1867, 12) 
		
	\end{quote}
	
	\noindent The new application of the PFL led Hankel to advocate a formal treatment of \textit{all} mathematics, which considered all operations as contentless, i.e., detached from any intuition of magnitudes, and operating solely within the domain of logical thought. In his last work, for example, he suggested an extension of this treatment to analysis:
	
	\begin{quote} \singlespacing
		In a similar way, just as for the system of numbers (or magnitudes) the whole numbers and their laws of operation were typical and declared to be permanent, so now for the system of functions the algebraic expressions -- those directly derived from the four fundamental operations -- must furnish the permanent laws. For just as the integers, in a certain sense, form the framework upon which all other kinds of numbers rest, more or less directly, so the algebraic functions, whose values can in fact be computed, are the typical forms for functions in general. (Hankel 1870, 101)
		
	\end{quote}
	
	\noindent Furthermore, Hankel also came to understand the PFL as a necessary condition for the legitimacy and usefulness of extending the system of functions:
	
	\begin{quote} \singlespacing
		Functions ... can become legitimate by subordinating themselves to a law. What, then, should this law be to which these functions must be subjected? It is not difficult to see that, if the concept of function is to be appropriate and meet the needs of the analyst, we must employ the general guiding principle that I have already established in my theory of complex numbers as the basis of all systematic progression in that field, and which I have called the principle of the permanence of formal laws. (1870, 101)
		
	\end{quote}
	
	\noindent On Hankel's view, new functions should be introduced only if they obey the laws of algebraic functions. More generally, the PFL should be understood, in a similar way, as a necessary condition for the legitimacy and usefulness of the doctrine of forms. That is, formal theories should be introduced only if they preserve the laws that hold in corresponding theories of magnitudes. For example, a theory of formal irrationals is legitimate and useful only if it satisfies the basic laws of the theory of actual irrationals. But Hankel thought that the PFL faced certain difficulties in this case, and likewise in the case of the hypercomplexes. Transcendental functions were considered problematic as well, but I will not have much to say about them in this paper. He noted that they are similar to the hypercomplexes in one respect: ``just as higher complex numbers can be conceived and constructed that do not subordinate themselves to the formal laws of arithmetic, so there can also be illegitimate functions which do not conform to the type of algebraic forms.'' (1870, 101) They are similar to the irrationals in another respect: ``the transcendental functions are not exhausted by a finite sequence of such operations [i.e., the four elementary operations], but require an infinite process. \ldots\ [they] correspond in a certain way to the irrational \ldots\ numbers'' (1870, 101sq). It seems fair, thus, to presume that the difficulties posed by transcendental functions might not go beyond those already pointed out in the other two cases. 
	
	In order to reconstruct Hankel's view about the irrationals, let us note that this rested on a particular criterion for scientific value: 
	
	\begin{quote} \singlespacing
		Every attempt to treat the irrational numbers formally and without the concept of magnitude must lead to the most highly abstruse and difficult artificialities, which, even if they can be carried out with complete rigor (which we have good reason to doubt), do not have a higher scientific value. For in general, it is the object of systematic science to become clear and conscious of the true foundations of the natural development of ideas, but not to desire to replace the organism, with its ever fresh power of production, with a dead and unproductive mechanism even if ingeniously constructed. (1867, 47)
		
	\end{quote}
	
	\noindent According to the criterion suggested here, in order to be valuable, a mathematical theory must not only be fully rigorous, it must be useful or productive, that is applicable to theories of magnitudes. The demand for productivity, in this sense, is supposed to balance epistemically the intuitively contentless operations of a formal theory: even though these operations are intuitively contentless, their results must be ultimately presented intuitively. Here Hankel encountered the following dilemma: the theory of irrational numbers is either not formal, if it requires operations with magnitudes presented in intuition, or not consistent, if it requires the unthinkable completion of an infinite series of operations. He put this as follows:
	
	\begin{quote}\singlespacing
		
		The operation required by an irrational number can be thought of as possible if the system of commensurable magnitudes, with which we operate, is a continuous one. Then there will be a series of number operations which, so long as it is not continued into the infinite, will always deviate from the one required, but comes closer to it without limit. The required operation can then be supposed to be carried out through the ideal completion of an infinite series of number operations. Now insofar as such an infinite continuation in its completeness is inconceivable and leads every time to a contradiction, we must claim that the irrational numbers are impossible, as long as no means has been found to present them in any other way than through division and multiplication and through number operations in general. (1867, 59)
	\end{quote}
	
	\noindent Since a theory cannot be productive if it is inconsistent, Hankel's verdict was that any valuable theory of the irrationals must be a theory of magnitudes: ``A systematic grasp of the irrational requires the concept of magnitude.'' (1867, 47) More particularly, he maintained that these numbers could only be the object of a theory of \textit{geometrical} magnitudes because its concept of continuity, intuitively presented through such magnitudes, makes the completion of an infinite series of operations entirely dispensable:
	
	\begin{quote} \singlespacing
		
		Geometry presents just such a means in its magnitude operations, which are independent of every concept of number, but only by looking at the concept of continuity, in which that very contradiction is hidden, as a given. Pure thought, detached from every intuition, cannot grasp the infinite; the formal theory of numbers cannot grasp the irrational. Intuition however requires the continuous: geometry proves the existence of the irrational. (1867, 59)
		
	\end{quote}
	
	\noindent Thus, according to Hankel, the difficulty that the irrationals raise for his doctrine of forms is that they are undetachable from intuitive content, and thus not formalizable. They are, in other words, \textit{essentially actual} numbers, and they can only be understood as expressions of relations between geometrical magnitudes. This is the reason why he also rejected the concept of limit as unable to support a formal theory of the irrationals: ``such a concept is based on the idea of small and large, which is entirely foreign to our development, and on the ordering of our numbers in a continuous series, which already involves the concept of extensive magnitude.'' (1867, 46) The dependence of the concept of limit on the concept of magnitude is, according to Hankel, ineliminable. Therefore, geometry gives the only permissible treatment of the irrationals, through an old and familiar application of the PFL: if a law holds for the actual rationals, then it shall hold for the actual irrationals as well. Or as he put it: ``exactly the same rules of operation apply to irrational numbers as to rational ones.'' (1867, 59) 
	
	Now, if a formal theory of irrationals cannot be reached through the new application of the PFL, the question arises why Hankel did not see this as a counterexample. Why did he not take the case of the irrationals to invalidate the PFL? Why did he not come to believe that such difficulties weaken its normative force? Why did he, instead, remain fully committed to it? The answer to these questions is to be found, I think, by carefully considering Hankel's view on the passage from complex to hypercomplex numbers. This will reveal the fact that, instead of simply restricting the scope of the PFL to avoid the problematic cases, which is what others like Schubert and Pringsheim would later do, Hankel adopted a conservative strategy that allows the revision of basic arithmetical laws, but retains sufficient normative force for the PFL.

	\section{Deliberative Conservatism}
	
	Hankel's view on the hypercomplexes is best understood on the background of Hamilton's own account. About a decade after the introduction of quaternions, Hamilton recalled one of his early attempts to extend the operation of multiplication of lines (or vectors) from two to three spatial dimensions, and noted that he did not find it satisfactory because ``it did not preserve the \textit{distributive principle} of multiplication'' (Hamilton 1853, 36). Other attempts were rejected for the same reason, as he found they all violated distributivity. Inconsistency with this law led him to abandon all these attempts, which indicates that Hamilton thought that if the usual laws of multiplication could be preserved, then they should be preserved:
	
	\begin{quote} \singlespacing
		With such preparations as I have described, I resumed
		(in	1843) the endeavour to adapt the general conception of triplets to the multiplication of lines in space, resolving to \textit{retain} the \textit{distributive} principle, with which some formerly conjectured systems had been inconsistent, and at first supposing that I \textit{could}
		preserve the \textit{commutative} principle also. (1853, 43)
		
	\end{quote}
	
	\noindent Supposing at first that multiplication obeys the law of commutativity, Hamilton realized that there is no satisfactory way to represent a product as a triplet. This gave him one reason to doubt that commutativity could be preserved. He then realized that multiplication would not yield unique results, which gave him an even stronger reason to doubt commutativity: 
	
	\begin{quote} \singlespacing
		
		The multiplication of \textit{lines} among themselves has been shewn to give \textit{different results}, according as the factors have been taken in one or in another \textit{order}; from which it follows, by still stronger reason, that the \textit{multiplication of quaternions} is \textit{not} generally a \textit{commutative} operation. (1853, 132)
		
	\end{quote}
	
	\noindent Nevertheless, instead of abandoning commutativity right away, Hamilton tried to preserve it. He eventually introduced the now well-known representation of a line as a quaternion $x+iy+jz+kw$, which takes the multiplication of quaternionic units to be non-commutative: $i^2 = j^2 = k^2 = -1$, $ij = -ji = k$, $jk = -kj = i$, and $ ki= - ik = j$. He then proved that commutativity fails for the product of any two ``rectangular vectors'' (i.e., mutually perpendicular lines): $\alpha\beta = - \beta\alpha$. Hamilton frankly acknowledged that this ``might at first seem strange'' (1853, 48) but noted that ``we are compelled, by considerations which appear more primary, to \textit{give up the commutative property} of multiplication, as not holding \textit{generally} for \textit{lines}.'' (1853, 51) Further, he showed that, unlike the law of commutativity, the laws of distributivity and associativity are preserved in quaternionic algebra, and moreover, he argued that they should both be preserved due to their usefulness in applications (1853, 494) Elsewhere, I have pointed out that Hamilton acknowledged not only Peacock's influence on his work, but also the more particular fact that in the passage from complex numbers to quaternions (and then to biquaternions), he actually applied Peacock's principle of the permanence of equivalent forms (Toader 2026). Hamilton believed that one should preserve the basic laws of arithmetic, but also that one should always weigh the reasons for their revision against the reasons for their preservation. This is deliberative conservatism.
	
	Accordingly, it is actually no surprise that, in a letter to his friend Mortimer O'Sullivan,  Hamilton described his own work on quaternions as ``at bottom quite conservative'':
	
	\begin{quote} \singlespacing
		
	You will I hope bear with me if I say, that it required a certain capital of scientific reputation, amassed in former years, to make it other than dangerously imprudent to hazard the publication of a work [i.e., Hamilton 1853] which has, although at bottom quite conservative, a highly revolutionary air. It was a part of the ordeal through which I had to pass, an episode in the battle of life, to know that even candid and friendly people secretly, or, as it might happen, openly, censured or ridiculed me, for what appeared to them my monstrous innovations. (quoted in Crowe 1967, 37).
	
\end{quote}

\noindent Hamilton's work was conservative, insofar as it maintained that all arithmetical laws should be preserved, but it was also revolutionary, insofar as it allowed that any of these laws could be revised. Having carefully studied Hamilton's writings, including the \textit{Lectures on Quaternions}, Hankel could not have missed the account given there of Hamilton's own deliberation on the properties of multiplication:
	
	\begin{quote}\singlespacing
		An interesting account of his efforts is contained in the preface to the Lectures. It shows that the main difficulty lay in discovering those properties of operations which, when abandoned, would least violate what I have called the Principle of the Permanence of Formal Laws of Arithmetic. Hamilton initially regarded the commutativity of multiplication as absolutely necessary, until numerous attempts convinced him that this property could not be maintained together with the distributive relation to the previously given addition of lines. He therefore resolved, in favor of the latter, to abandon the commutative property of multiplication. (Hankel 1867, 195; also 105)
	\end{quote}
	
	\noindent Hankel's emphasis on a comparative evaluation of hypercomplex systems strongly suggests an interpretation of the PFL in line with deliberative conservatism. On his view, Hamilton's task was to preserve as many of the laws of arithmetic as possible, while considering the consequences of abandoning any of them. Closely following Hamilton, Hankel acknowledged that the application of the PFL admits of exceptions:
	
	\begin{quote} \singlespacing
		
		This principle will guide our steps in what follows; but one is not allowed to apply it without exception and universally; we shall allow universal application only to the definition of necessary and sufficient rules, to the extent that these are independent of one another; still we will not allow it to restrict us too severely, particularly, we shall not unconditionally assume the commutativity of our operations -- since it has proved to be a scientific necessity to consider operations which conform to all the rules of arithmetic multiplication, with this single exception. (1867, 11sq)
		
	\end{quote}
	
	\noindent Hankel added that through the application of the PFL, the laws of addition and multiplication governing hypercomplex systems are ``modeled after the known laws of actual addition and multiplication in the doctrine of magnitudes \ldots\ with some freedom (with regard to commutativity)'' (1867, 33). If the commutativity is omitted, ``the equation $xx = -1$ can have more roots than just $\pm i$, as it indeed does in the theory of quaternions'' (1867, 70). However, just how far to carry the PFL is always a matter for deliberation: ``in order to make the laws of arithmetic multiplication permanent \textit{as far as possible}, one could, e.g., preserve commutativity, and stipulate $i_n \cdot i_p = i_p \cdot i_n$'' (1867, 102, my emphasis), in which case, as Hamilton had already shown, distributivity would have to be omitted. Alternatively, as in Hamilton's own system, ``by abandoning that permanence [of commutativity, one could] stipulate $i_n \cdot i_p = - i_p \cdot i_n$'' (1867, 102) thereby preserving distributivity instead. There is, indeed, textual evidence that Hankel took seriously the possibility of a commutative hypercomplex system:
	
	\begin{quote}\singlespacing
		
		We have in our text only thought of such complex numbers whose multiplication satisfies the associative and distributive principle. In and of itself, however, nothing stands opposed to letting one of these laws fall partially and instead, for example, preserving the commutative one. In fact, H. Scheffler ... has set up a type of complex numbers for treating spatial relations, which are entirely commutative in their multiplication, but not subject to the distributive law in general. (1867, 105)
	\end{quote}
	
	\noindent The very fact that Hankel understood Hermann Scheffler's algebraic calculus of directed quantities (\textit{gerichtete Größe}), which had been assumed commutative but not distributive (Scheffler 1846, Scheffler 1851), as the basis for a possible alternative to Hamilton's own system, supports my interpretation of the PFL as an expression of deliberative conservatism, for it shows that Hankel regarded the choice between commutativity and distributivity as a genuine deliberative trade-off. 
	
	Now, if this is the right interpretation, then it can also explain why Hankel did not see the irrationals as a counterexample to the PFL. The explanation is as follows: although the principle stipulates that the laws that hold for the actual irrationals should also hold for the formal irrationals, it allows that they can be abandoned if the reasons for preserving them are outweighed by the reasons for abandoning them. And clearly enough, Hankel thought that those laws must be abandoned, for the inconsistency that he saw in the completion, within the domain of logical thought, of an infinite series of operations makes the preservation of the laws logically impossible. Thus, the difficulty raised by this case is not a violation that invalidates the PFL, which explains why he remained strongly committed to his principle. The doctrine of forms is not only worth developing, but can actually be developed in the way he thought should be developed, through his new and independent application of the PFL.
	
	Furthermore, it should be clear why the sort of deliberation licensed by the PFL underlies an important part of Hankel's anti-Kantian program, though a part that seems to have been neglected so far. My interpretation implies that, on his view, the basic laws of arithmetic are not unrevisable; rather, they are subject to revision: each or even all of them can be omitted or abandoned, if there are strong reasons for doing so. In this sense, Hankel's conservatism is definitely anti-Kantian, for it is directly opposed to the kind of conservatism associated to a Kantian view of mathematics, according to which arithmetic is epistemically privileged, unconditionally valid, and immune to revision (see Callanan 2019, Kjosavik 2022). Kant himself had treated arithmetic as a theory of magnitudes that was the foundation of his conception of mathematical cognition (Kant 1781; see Parsons 1969, Sutherland 2021). As is well known, his view was further elaborated by his contemporary followers, like Johann Schultz, who explicitly argued, e.g., that the commutativity and associativity of addition are ``axioms'' of arithmetic that are indispensable (``\textit{unentbehrlich}'') to any extension of the theory, such as ``higher arithmetic'', ``general'' or ``pure'' \textit{Mathesis}, i.e., algebra (Schultz 1791, 225; see Martin 1972, Parsons 2012, von Plato 2016). This view, or more precisely one that can be interpreted along these Kantian lines, is detectable in the works of later mathematicians, like Schubert and Pringsheim, who turned against Hankel.

	\section{Against Hankel}

	A founding member of the \textit{Deutsche Mathematiker-Vereinigung} and Professor at the Johanneum Gymnasium in Hamburg, Hermann Schubert is best known today, among mathematicians and physicists, for his theory of enumerative geometry, which famously became the subject of Hilbert's 15th problem (Schubert 1879; see Bethea and Brazelton 2025). In several of his more popular writings about mathematics, Schubert expressed the following philosophical view: 
	
	\begin{quote}\singlespacing
		The intrinsic character of mathematical research and knowledge is based essentially on three properties: first, on its conservative attitude towards the old truths and discoveries of mathematics; secondly, on its progressive mode of development, due to the incessant acquisition of new knowledge on the basis of old knowledge; and thirdly, on its self-sufficiency and its consequent absolute independence. (Schubert 1896, 294). 
	\end{quote}
	
	\noindent The kind of conservatism that he defended was emphatically anti-revisionary. This was especially apparent in his conception of arithmetic, which he said was developed ``wholly and alone in virtue of the logical application of the monistic principle of arithmetic, the principle of no exception'' (Schubert 1894, 579). He formulated this principle in his contribution to the first volume of the \textit{Encyklopädie der mathematischen Wissenschaften}, edited by Felix Klein, but not before recalling in a footnote the two precursors of the principle: 
	
	\begin{quote}\singlespacing
		
		The principle of \textit{permanence}, which in the text here is given the form suitable for the extension of the concept of number, was first expressed in its most general form by H. Hankel [...] after G. Peacock had already emphasized the necessity of a purely formal mathematics and, in connection with this, a principle from which the principle of permanence proceeds by extension. (Schubert 1898, 11)
		
	\end{quote}
	
	\noindent Elsewhere, I have discussed Schubert's four-part formulation of his principle, and Frege's and Peano's critical remarks on some of the four tasks it describes (Toader 2021). Here I want to briefly emphasize that the principle of no exception, or as Schubert also called it, the principle of exceptionlessness (\textit{Ausnahmslosigkeit}), implies a rejection of Hankel's deliberative conservatism. As its name already suggests, Schubert's principle denies outright the possibility of revising the basic laws of arithmetic. As a consequence, it cannot sit well with typical extensions of number systems beyond the complexes. This led Schubert to reject Hamilton's quaternions as genuine numbers: 
	
	\begin{quote}\singlespacing
		With respect to quaternions, which many might be disposed to regard as new numbers, it will be evident that though quaternions are valuable means of investigation in geometry and mechanics they are not numbers of arithmetic, because the rules of arithmetic are not unconditionally applicable to them. (1894, 578)
	\end{quote}
	
	\noindent This rejection of quaternions makes the obvious assumption that something can be a number only if it unconditionally obeys all the laws of arithmetic -- an assumption Schubert evidently found less problematic than he found Hankel's deliberative conservatism. The latter, as we have seen, treats the extension of number systems as genuinely open, in the following sense: although it takes the preservation of laws as a requirement, exactly what laws are to be preserved is always subject to rational deliberation. But Schubert was not alone in criticizing Hankel. 
	
	Alfred Pringsheim, Professor at the University of Munich and President of the \textit{Deutsche Mathematiker-Vereinigung}, who also contributed articles to Klein's \textit{Encyklopädie} on irrational numbers and on infinite processes (in 1898, and again in 1904), criticized Hankel as well. In his widely read \textit{Vorlesungen über Zahlen- und Funktionenlehre} (5 volumes, published between 1916--1932), Pringsheim characterized Hankel's attempt to give a formal treatment to the irrationals as inadequate, and rejected the conclusion that this is logically impossible: ``it is extremely remarkable that even the creator of a \textit{purely formal theory of rational numbers} has shown so little understanding for the corresponding further development of the number concept'' (Pringsheim 1898, 57). The irrationals, treated formally, present no special obstacles. Of course, unlike Hankel, Pringsheim was in a position to cite works that had provided rigorous definitions of the irrationals (e.g., by Cantor, Dedekind, and others).
	
	Indeed, Pringsheim endorsed the first part of Hankel's anti-Kantian program: the formalization of all mathematics, its detachment from any intuition of magnitudes. As such, he rejected Du Bois-Reymond's critique of the detachment of irrationals from geometrical magnitudes. However, Pringsheim was not so adamant about the second part of that program: Hankel's deliberative conservatism. He reformulated the PFL as a \textit{transfer principle} (\textit{Übertragungsprinzip}), explicitly distancing himself from both Hankel's terminology and his understanding of it:
	
	\begin{quote}\singlespacing
		To establish how these new numbers are to be \textit{ordered} or, alternatively, \textit{incorporated} into the ordered sequence of already existing numbers and how we are to \textit{calculate} with them, we shall make use of the \textit{transfer principle} which (following Hankel) is usually (but not very felicitously) called the principle of `\textit{permanence}', and we shall make use of it in what I regard as a notably improved form which bestows on it the character of a certain logical necessity. (Pringsheim 1916, VIIsq)
	\end{quote}
	
	Not only was ``permanence'' a misnomer, on his view, but Pringsheim believed that the PFL lacked a kind of logical necessity. This is the main reason why he took his own transfer principle to constitute a notable improvement. But since, as we have seen in the previous section, Hankel himself had come to think of the PFL as a necessary (though surely not sufficient) condition for the legitimacy and usefulness of formal theories, then we should ask in what sense Pringsheim's view differs from Hankel's. 
	
	In order to explain Pringsheim's view, let us note two elements of the transfer principle. The first is what he called transcription (\textit{Umschreibung}): when new numbers are introduced, new number signs are introduced, and this must be done in such a way that some of the new signs replace the old notation for the already existing numbers and the rules governing the latter can be transcribed into the new notation. In his own words:
	
	\begin{quote}\singlespacing
		We shall introduce \textit{new number signs}, but only to such an extent that a subset of them represents signs for \textit{already existing} numbers. The latter are therefore already governed by certain rules establishing their succession and defining the arithmetical operations for them, and these rules can without further ado be \textit{trans}cribed into the new notation. (1916, VIIsq)
	\end{quote}
	
	The second element of the principle is what he called extension (\textit{Ausdehnung}): the transcribed rules, which initially cover only a subset of the new signs, are extended to cover all of them. These two elements, the transcription and extension of rules, together with a proof that the extension is consistent, confer the transfer principle the kind of logical necessity that Pringsheim did not see in the PFL, and by which he seems to have meant epistemic necessity, or in any case, the character of a necessary condition for the possibility of the consistent manipulation of signs:
	
	\begin{quote}\singlespacing
		If we are not to allow complete confusion in the manipulation of the total supply of newly created signs, we have hardly any choice but to extend the rules already governing part of it to the totality by definition, and to legitimize this step by proving that the stipulations we have made satisfy the requirements to be met by them without contradiction. (1916, VIIsq)
	\end{quote}
	
	This can be understood as a transcendental argument for the transfer principle. If so, then despite his endorsement of the first part of Hankel's anti-Kantian program --- formalization of all mathematics --- Pringsheim turns out to be closer to Schubert's conservatism: unless all laws of arithmetic are preserved, one cannot prevent the manipulation of signs from being completely confused or incoherent. The difference here is significant: the PFL, as a necessary condition for legitimizing formalization, leaves room for deliberation; the transfer principle, as a necessary condition for coherent sign manipulation, is to determine any extension uniquely.
	
	If Pringsheim's transfer principle is understood in this way, then its application to the irrationals raises no problems. The procedure used to expand the natural numbers to the rationals serves perfectly well the introduction of irrational numbers: new signs are introduced, the rules for the rationals are transcribed into the new notation, extended, and proved consistent. ``By simply \textit{transferring} the laws that exist for the ordering and combination of the rational numbers \textit{into the new notation} one obtains without further ado a corresponding complex of \textit{already secured} rules for this new type of number signs'' (1916, VIII). Pringsheim applied the same procedure to the complex numbers, although this requires an ``appropriate modification'' (Pringsheim 1921, V) which is worth noting.
	
	Deviating from the standard textbook approach of defining complex numbers as pairs of reals, his own account starts by assigning to each real number $\alpha$ a new imaginary number $\bar{\alpha}$, with the defining stipulation: $\bar{\alpha}\cdot\bar{\beta} = -(\alpha\beta)$. The meaning of the sign combination $\bar{\alpha}\cdot\bar{\beta}$ is given by the already existing real number $-(\alpha\beta)$, and all the laws of the reals are transferred as definitions of operations on the imaginaries. The set of all real and imaginary numbers is then given an order, but Pringsheim presented this only as a formal possibility designed to preserve the existing order relations: for any arbitrarily small positive real number $\epsilon$ and any arbitrarily large positive real $\omega$,
	\[
	-\omega < -\epsilon < \overline{(-\omega)} < \overline{(-\epsilon)} < 0 < \bar{\epsilon} < \bar{\omega} < \epsilon < \omega,
	\]
	and the sum of two imaginary numbers is defined by $\bar{\alpha}+\bar{\beta}=\overline{\alpha+\beta}$. Since the sum $\alpha+\bar{\beta}$ cannot be classified as either real or imaginary under the earlier rules, Pringsheim argued that such sign combinations must be introduced as complex numbers, of which the reals and imaginaries become special cases, and whose laws are then transferred to the complexes, he emphasized, ``with compelling necessity'' (1921, 526).
	
	The transfer principle is applied in a similar manner at each step in the construction of the number system. However, Pringsheim never applied it, and does not seem to have ever been interested in applying it, to extensions beyond the complexes. In all likelihood, following Schubert, he was content to simply relinquish hypercomplex numbers. Surely one can imagine someone like Hankel complaining about this with some irony: it is extremely remarkable that ``the most eminent representative of the arithmetical school'' (Hahn 1919, 66) could show so little understanding for the further development of the number concept... Quite so, for if the transfer principle has the kind of necessity that Pringsheim wanted -- as a necessary condition for the consistent manipulation of any sign system -- then it is not clear why it should not apply to quaternions as well. One explanation for this is that he understood his own transfer principle in a way that is closer to Schubert's anti-revisionary conservatism than to Hankel's deliberative conservatism.

	\section{Conclusion}
	
	This paper has analyzed Hankel's applications of the principle of the permanence of formal laws (PFL) and argued that this constitutes a two-part anti-Kantian program in the philosophy of mathematics. The first part, presented in Section 2, is based on Hankel's new and independent application of the PFL for guiding and legitimizing formal theories: if a law holds for actual mathematical objects, it should hold for their formal counterparts as well. This application, which underlies Hankel's doctrine of forms, advocates a formal treatment of all of mathematics, detached from any intuition of magnitudes -- a broadly formalist stance that, as I noted at the outset, was shared by others. The second part, discussed in Section 3, concerns deliberative conservatism -- a conservative strategy that Hankel inherited from Peacock and Hamilton, which permits the revision of laws when the reasons for revision outweigh the reasons for preservation. This strategy explains why Hankel remained committed to the PFL despite the difficulties raised by the irrationals and hypercomplexes: these are not counterexamples that invalidate the principle, but cases in which the deliberation it licenses yields the conclusion that some or even all laws cannot be preserved. Section 4 then examined the criticisms of Schubert and Pringsheim, who reinterpreted the PFL to align it with their own anti-revisionary conservatism. 
	
	But the philosophical wedge between these two kinds of conservatism runs deeper than might be suggested by the differences between Hamilton and Hankel, on the one hand, and Schubert and Pringsheim, on the other. Elsewhere, I have argued that deliberative conservatism is philosophically grounded in Hume's conception of the laws of reasoning, as articulated in section 1.4.4.1 of his \textit{Treatise} (Hume 1739, see Toader 2026). These laws, according to Hume, are permanent, but in principle revisable when sufficiently strong reasons are brought against them. On his view, as I understand it, the laws of probable (or inductive), as well as the laws of demonstrative (or deductive) reasoning must be preserved to the furthest extent possible. This allows for exceptions to these laws, for as Hume noted, although their violation is practically unsustainable in general, it is in principle not impossible to violate them. Anti-revisionary conservatism, of the sort defended by Schubert and Pringsheim, treats the basic laws of arithmetic as unconditionally binding, immune to revision by any mathematical consideration whatsoever. This view is best understood on the background of Kant's philosophy of mathematics, which treated arithmetic not as open to revision, but as furnishing necessary and unrevisable conditions for the possibility of mathematical knowledge. As I noted in the paper, this was endorsed and further articulated by some of Kant's contemporaries, especially Schultz, who considered the axioms of arithmetic to be not only unconditionally valid, but unrevisable. Later, the Marburg neo-Kantians transformed Kant's conception of arithmetic, and did so in ways that deviate from Schultz. Already Cohen, for example, treated the commutativity of addition as a theorem, rather than an axiom (Cohen 1902, 161). Since it is derived from deeper transcendental laws, this arithmetical theorem is conditional on these laws and thus could in principle be revised. A closer analysis of this view must, however, be deferred to future work.
	
	Some questions about the relationship between the two parts of Hankel's anti-Kantian program, and more generally between formalization and conservatism, also deserve further investigation. For instance, suppose that one is interested in the application of the PFL to formal systems, rather than to theories of magnitudes, so that the principle stipulates that if a formal system has a certain property (consistency, completeness, categoricity, decidability, etc.), then its formal extensions should have it as well. Which kind of conservatism better supports this specifically formalist application? On the one hand, Kantian conservatism seems to offer a more stable foundation for formalization, but the rigidity it demands may preclude precisely the extensions that Hankel's doctrine of forms was designed to accommodate. Deliberative conservatism, on the other hand, is more flexible and hospitable to formal extensions, but may undermine the uniqueness that one might want them to possess. In fact, it is not clear that even Pringsheim's insistence on the logical necessity of the transfer principle is enough to determine a unique extension of any given sign structure (Hahn 1919; see Toader 2021). But in addition, one should ask if either of these kinds of conservatism was explicitly or implicitly adopted by later formalists. Hilbert's formalism, for instance, insisted on consistency proofs as the ultimate criterion of legitimacy, but his readiness to entertain new axiom systems seems to indicate a deliberative inclination. He suggested that new axioms could be justified by what he called ``systematic advantages (principle of the permanence of laws, further possibilities of construction, etc.).'' (Hilbert 1929) Should this be understood as an expression of deliberative conservatism, or as a commitment to Kantian conservatism? This question points toward an investigation that deserves another paper.
	
\section{Acknowledgments}

The Austrian Science Fund (FWF) is gratefully acknowledged for financial support through 
Grant 10.55776/PAT3440123.

	\section{References}
	
	\begin{description}
		\singlespacing
		\setlength\itemsep{0.4em}

\item Bethea, C. and Brazelton, T. 2025. The evolution of enumerative geometry: a narrative from classical problems to enriched invariants.  https://doi.org/10.48550/arXiv.2510.04275

\item Bolzano, B. 1810. \textit{Beyträge zu einer begründeteren Darstellung der Mathematik}. Prague: Caspar Widtmann. Eng. tr. in: \textit{The Mathematical Works of Bernard Bolzano}. Oxford, 2004

\item Callanan, J. 2019. Methodological conservativism in Kant and Strawson. \textit{British Journal for the History of Philosophy} 27: 422–442

\item Cohen, H. 1902. \textit{Logik der reinen Erkenntniss}. Berlin: Bruno Cassirer. 3rd ed. 1922

\item Crowe, M. J. 1967. \textit{A History of Vector Analysis: The Evolution of the Idea of a Vectorial System}. Notre Dame: University of Notre Dame Press. 2nd ed. 1985, Dover Publications

\item Förstemann, W. A. 1817. \textit{Über den Gegensatz positiver und negativer Größen}. Nordhausen: Happach. Reprint, Wiesbaden: LTR-Verlag, 1971

\item Frege, G. 1884. \textit{Die Grundlagen der Arithmetik}. Breslau: Koebner. Eng. tr.: \textit{The Foundations of Arithmetic}. Oxford: Blackwell, 1950

\item Frisch, U. and Villone, B. 2014. Cauchy's almost forgotten Lagrangian formulation of the Euler equation for 3D incompressible flow. \textit{European Physical Journal H} 39: 325–351

\item Hahn, H. 1919. Review of A. Pringsheim, \textit{Vorlesungen über Zahlen- und Funktionenlehre}. \textit{Gött. Gel. Anz.} 9: 321–347. Reprinted in: \textit{Empiricism, Logic, and Mathematics}. Dordrecht: Reidel, 1980, 51–72

\item Hamilton, W. R. 1853. \textit{Lectures on Quaternions}. Dublin: Hodges and Smith.

\item Hankel, H. 1861. \textit{Zur allgemeinen Theorie der Bewegung der Flüssigkeiten}. Göttingen: Dieterich

\item Hankel, H. 1867. \textit{Theorie der complexen Zahlensysteme}. Leipzig: Voss

\item Hankel, H. 1869. \textit{Zur Geschichte der Mathematik in Alterthum und Mittelalter}. Leipzig: Teubner

\item Hankel, H. 1870. \textit{Untersuchungen über die unendlich oft oscillirenden und unstetigen Functionen}. Tübingen: Fues

\item Hankel, H. 1874. \textit{Die Entwicklung der Mathematik in den letzten Jahrhunderten}. Tübingen: Fues

\item Hilbert, D. 1929. Probleme der Grundlegung der Mathematik. \textit{Mathematische Annalen} 102: 1–9. Eng. tr. in: Mancosu, P. (ed.) \textit{From Brouwer to Hilbert}. Oxford, 1998, 227–233

\item Hume, D. 1739. \textit{A Treatise of Human Nature}, ed. 2007, Clarendon Press

\item Kant, I. 1781. \textit{Kritik der reinen Vernunft}. Riga: Hartknoch. Eng. tr.: \textit{Critique of Pure Reason}. Cambridge: Cambridge University Press, 1998

\item Kjosavik, F. 2022. Kant on the possibilities of mathematics and the scope and limits of logic. \textit{Inquiry} 65: 683–706

\item Martin, G. 1972. \textit{Arithmetic and Combinatorics: Kant and His Contemporaries}. Carbondale: Southern Illinois University Press, 1985

\item Parsons, C. 1969. Kant's philosophy of arithmetic. In: \textit{Mathematics in Philosophy}. Ithaca: Cornell University Press, 1983

\item Parsons, C. 2012. \textit{From Kant to Husserl}. Cambridge, MA: Harvard University Press

\item Peacock, G. 1830. \textit{A Treatise on Algebra}. Cambridge: J. \& J. J. Deighton

\item Peacock, G. 1833. Report on the recent progress and present state of certain branches of analysis. \textit{Report of the Third Meeting of the British Association for the Advancement of Science} 3: 185–352

\item Peacock, G. 1842. \textit{A Treatise on Algebra}, Vol. I, 2nd ed. Cambridge: J. \& J. J. Deighton

\item Peacock, G. 1845. \textit{A Treatise on Algebra}, Vol. II, 2nd ed. Cambridge: J. \& J. J. Deighton

\item Pringsheim, A. 1898. Irrationalzahlen und Konvergenz unendlicher Prozesse. In: \textit{Encyklopädie der mathematischen Wissenschaften mit Einschluss ihrer Anwendungen}, Vol. I. Leipzig: Teubner, 47–146

\item Pringsheim, A. 1916. \textit{Vorlesungen über Zahlen- und Funktionenlehre}, Vol. 1, Part 1. Leipzig: Teubner

\item Pringsheim, A. 1921. \textit{Vorlesungen über Zahlen- und Funktionenlehre}, Vol. 1, Part 2. Leipzig: Teubner

\item Salazar, M. S. 2023. \textit{The Development of Conceptualizations of Number Domains as the Basis for the New Rigor in 19th-Century Mathematics}. PhD thesis, Federal University of Rio de Janeiro

\item Scheffler, H. 1846. \textit{Über das Verhältniss der Arithmetik zur Geometrie, insbesondere über die geometrische Bedeutung der imaginären Zahlen}. Braunschweig: Leibrock

\item Scheffler, H. 1851. \textit{Der Situationscalcul}. Braunschweig: Leibrock

\item Schubert, H. 1879. \textit{Kalkül der abzählenden Geometrie}. Leipzig: Teubner

\item Schubert, H. 1894. Monism in arithmetic. \textit{The Monist} 4: 561–579

\item Schubert, H. 1896. On the nature of mathematical knowledge. \textit{The Monist} 6: 294–305

\item Schubert, H. 1898. Grundlagen der Arithmetik. In: \textit{Encyklopädie der mathematischen Wissenschaften mit Einschluss ihrer Anwendungen}. Leipzig: Teubner, 1–27

\item Schubring, G. 2005. \textit{Conflicts Between Generalization, Rigor, and Intuition}. New York: Springer

\item Schultz, J. 1791. \textit{Prüfung der kantischen Critik der reinen Vernunft}, Vol. 1. Frankfurt and Leipzig

\item Sutherland, D. 2021. \textit{Kant's Mathematical World: Mathematics, Cognition, and Experience}. Cambridge: Cambridge University Press

\item Toader, I. D. 2021. Permanence as a principle of practice. \textit{Historia Mathematica} 54: 77–94

\item Toader, I. D. 2026. Peacock's principle as a conservative strategy. \textit{Archive for History of Exact Sciences} 80: 9

\item von Plato, J. 2016. In search of the roots of formal computation. In: Gadducci, F. and Tavosanis, M. (eds.) \textit{HaPoC 2015}, 300–320
		
	\end{description}
	
\end{document}